\documentclass[ECP]{ejpecp} 

\SHORTTITLE{Random matrix products}

\TITLE{Nonconventional random matrix products
\thanks{S.Sodin was supported in part by the European Research Council starting grant SPECTRUM (639305) and by a Royal Society Wolfson Research Merit Award.}}

\AUTHORS{%
  Yuri~Kifer\footnote{Institute of Mathematics, Hebrew University of Jerusalem.
    \EMAIL{kifer@math.huji.ac.il}}
  \and 
  Sasha~Sodin\footnote{School of Mathematical Sciences, Queen Mary University of
  London \& School of Mathematical Sciences, Tel Aviv University. \EMAIL{asodin@qmul.ac.uk}}}


\AMSSUBJ{60B20; 60F15; 60F10; 82B44} 

\SUBMITTED{March 25, 2018} 
\ACCEPTED{May 29, 2018} 

\ARXIVID{1803.09221} 

\VOLUME{0}
\YEAR{2018}
\PAPERNUM{0}
\DOI{10.1214/YY-TN}

\ABSTRACT{Let $\xi_1,\xi_2,...$ be independent identically distributed random variables
and $F:\bbR^\ell\to SL_d(\bbR)$ be a Borel measurable matrix-valued function.
Set $X_n=F(\xi_{q_1(n)},\xi_{q_2(n)},\\...,\xi_{q_\ell(n)})$ where
$0\leq q_1<q_2<...<q_\ell$ are increasing functions taking on integer values
on integers. We study the asymptotic behavior as $N\to\infty$ of
the singular values of the random matrix product $\Pi_N=X_N\cdots X_2X_1$ and
show, in particular, that (under certain conditions) $\frac 1N\log\|\Pi_N\|$
converges with probability one as $N\to\infty$. We also obtain  similar results
for such products when $\xi_i$ form a Markov chain. The essential difference
from the usual setting appears since the sequence $(X_n,\, n\geq 1)$ is
long-range dependent and nonstationary.}

\newcommand{\cE}{{\mathcal E}}
\newcommand{\cF}{{\mathcal F}}

\newcommand{\Om}{{\Omega}}
\newcommand{\om}{{\omega}}
\newcommand{\ve}{{\varepsilon}}
\newcommand{\del}{{\delta}}
\newcommand{\Del}{{\Delta}}
\newcommand{\gam}{{\gamma}}
\newcommand{\Gam}{{\Gamma}}
\newcommand{\vf}{{\varphi}}

\newcommand{\sig}{{\sigma}}
\newcommand{\al}{{\alpha}}

\newcommand{\ka}{{\kappa}}
\newcommand{\la}{{\lambda}}

\newcommand{\bbR}{{\mathbb R}}

\newcommand{\bbI}{{\mathbb I}}

\begin{document}

\section{Introduction}\label{sec1}\setcounter{equation}{0}

Products $\Pi_N=X_N\cdots X_2X_1$ of
random matrices $X_1,X_2,...$ are extensively studied for more than half a
century now. In the pioneering work \cite{FK}, it was shown that when $X_1,X_2,...$ form a stationary
sequence with $E\ln^+\| X_1\|<\infty$
then the limit $\gam_1=\lim_{N\to\infty}\frac 1N\ln\|\Pi_N\|$ exists with probability
one. Later, the more general Kingman's subadditive ergodic theorem became available and it
yielded the above result as a corollary. Applying it to actions on the exterior products,  the
result was extended to all the singular values of $\Pi_N$, thus leading to the Oseledets
multiplicative ergodic theorem.

In this paper we study similar questions for products of certain nonstationary sequences
of random matrices. Namely, we start with a sequence of i.i.d.\ random variables $\xi_1,\xi_2,
...$ and a Borel measurable matrix valued function $F:\bbR^\ell\to SL_d(\bbR)$ along
with integer valued functions $0\leq q_1<q_2<...<q_\ell$, and form the random matrices
$X_n=F(\xi_{q_1(n)},\xi_{q_2(n)},...,\xi_{q_\ell(n)})$. In particular, we allow arithmetic progressions $q_i(n)=in,\, i=1,...,\ell$.
The sequence $X_1,X_2,...$ is long range dependent and is not stationary, and so the study
of the asymptotic behavior
as $N\to\infty$ of the product $\Pi_N=X_N\cdots X_2X_1$ is not described by the standard
results mentioned above.  Still, we show that $\lim_{N\to\infty}\frac 1N\ln
\|\Pi_N\|$ exists with probability one and applying this to exterior products we will
obtain corresponding results for all the singular values of $\Pi_N$. Similar results are
obtained also for such products when $X_n=F(\xi_{n},\xi_{2n},...,\xi_{\ell n})$
and $\xi_i$ form a Markov chain satisfying certain conditions of the
type of uniform geometric ergodicity.

\medskip
The motivation for this paper is twofold. On one hand, it comes from the vast body of
research on products of random matrices mentioned above (see \cite{BL} and \cite{BQ}). In particular,  our results provide a non-trivial family of random discrete Schr\"odinger equations $\psi_{n+1}=(\la-V_n)\psi_n-\psi_{n-1}$ which are not metrically transitive and yet the asymptotics of solutions can be described, where, as usual,
$\Del\psi(n)=-(\psi(n+1)+ \psi(n-1))$ is viewed as a discrete counterpart of the Laplacian.
  In our case, $V_n=\vf(\xi_{q_1(n)},
\xi_{q_2(n)},...,\xi_{q_\ell(n)})$ and $\xi_1,\xi_2,...$ are, say, i.i.d.\ random variables.

On the other hand, our motivation stems from the series of papers, originating in Furstenberg's proof of the Szemer\'edi theorem, on nonconventional
ergodic and limit theorems which dealt with the sums of the form $\sum_{n=1}^N
\vf(\xi_{q_1(n)},\xi_{q_2(n)},...,\xi_{q_\ell(n)})$ (see, for instance, \cite{HK} and
references therein). Our results can be viewed as a counterpart of the nonconventional strong
law of large numbers in the multiplicative setting.

\section{Preliminaries and main results}\label{sec2}\setcounter{equation}{0}

\subsection{I.i.d.\ case} Let $\xi_1,\xi_2,...$ be i.i.d.\ random variables, and let $F:\bbR^\ell\to SL_d(\bbR)$ be a Borel measurable matrix valued function where $\ell>1$ (since for $\ell=1$ the results of this paper are well known). Our setup also includes
an $\ell$-tuple of strictly increasing nonnegative functions $q_1<q_2<...<q_\ell$ taking on integer
 values on integers with $q_1(1)\geq 1$. Set $X_n=F(\xi_{q_1(n)},\xi_{q_2(n)},...,
 \xi_{q_\ell(n)})$ and observe that each $X_n,\, n\geq 1$ has the same distribution,
 since each $\ell$-tuple $\xi_{q_1(n)},\xi_{q_2(n)},..., \xi_{q_\ell(n)}$
 has the same distribution as $\xi_1,\xi_2,...,\xi_\ell$. Denote by $\mu$ the distribution of $X_1$ and by $G_\mu$ the support of $\mu$. We will need the following

\begin{assumption}\label{ass2.1} \hfill
\begin{enumerate}
\item[(i)] $G_\mu$ is strongly irreducible, i.e.\ there does not exist a finite union of proper subspaces of $\mathbb R^d$ that is preserved as a set by all matrices from $G_\mu$ (see \cite{BL}).
\item[(ii)] For some $\al>0$,
\begin{equation}\label{2.1}
E \|X_1\|^\alpha < \infty.
\end{equation}
\item[(iii)]for any $\sigma>0$ there exists $n_0(\sigma)$ such that for all $n\geq n_0(\sigma)$,
\begin{equation}\label{2.2}
 q_{i+1}(n)\geq q_i(n+[\sigma\ln n]),\quad i=1,...,\ell-1.
\end{equation}
\end{enumerate}
\end{assumption}
Clearly, (\ref{2.2}) is satisfied, for instance, in the arithmetic progression case
$q_i(n)=in,\, i=1,...,\ell$.

Recall that the singular values $s_1(g)\geq s_2(g)\geq...
\geq s_d(g)\geq0$ of a $d \times d$ matrix $g$ are the square roots of the eigenvalues $s^2_i(g)$ of $g^*g$. The first singular value $s_1(g)$ is the Euclidean operator
norm of $g$,
\[
s_1(g)=\max_{x\in\bbR^d\setminus\{0\}}\frac {\| gx\|}{\| x\|}=\| g\|.
\]
If $X\in SL_d(\bbR)$ then $1=s_1(X)s_2(X)\cdots s_d(X)\leq s_1^{d-1}(X)s_d(X)$, and so $\| X^{-1}\|
=s_d^{-1}(X)\leq s_1^{d-1}(X)=\| X\|^{d-1}$. Hence, (\ref{2.1}) implies also that
\[
E\| X_1^{-1}\|^{\al'}<\infty\quad\mbox{with}\quad\al'=\frac \al{d-1},\, d>1.
\]
Since $F\equiv 1$ if $d=1$ and the problems discussed here become trivial then, we assume without loss of generality that $d>1$.

Let $Y_1,Y_2,...$ be an i.i.d.\ sequence of random
matrices having the distribution $\mu$, and so satisfying (i) and (ii) of Assumption
\ref{ass2.1} with $Y_1$ in place of $X_1$. Hence (cf.\ \cite{BL, BQ}),
 the limits
 \begin{equation}\label{2.3}
 \gam_i=\lim_{N\to\infty}\frac 1N\ln s_i(Y_N\cdots Y_2Y_1),\, i=1,...,d,
\end{equation}
 exist with probability one; in particular, $\gam_1=\lim_{N\to\infty}\frac 1N\ln\| Y_N\cdots Y_2Y_1\|$. The following theorem asserts that the similar result holds
 true for $\Pi_N=X_N\cdots X_2X_1$.
 \begin{theorem}\label{thm2.2}
 Suppose that Assumption \ref{ass2.1} holds true. Then with probability one
 \begin{equation}\label{2.4}
 \lim\limits_{N\to\infty}\frac 1N\ln s_i(\Pi_N)=\gam_i,\, i=1,...,d,
 \end{equation}
 where $\gam_1,...,\gam_d$ are the same as in (\ref{2.3}). In particular, $\lim\limits_{N\to\infty}\frac 1N\ln\|\Pi_N\|=\gam_1$.
 \end{theorem}

 Observe that it suffices to prove Theorem \ref{thm2.2} only for the largest singular value,
 i.e.\ for $i=1$. Indeed, observing that (i) and (ii) of Assumption \ref{ass2.1} remain valid
 for the exterior powers $\wedge^i\Pi_N,\, i=1,...,d$ of $\Pi_N$ (defined by $\wedge^i\Pi_N(x_1\wedge...
 \wedge x_i)=\Pi_Nx_1\wedge...\wedge\Pi_Nx_i$) if this were true for
 $\Pi_N$ itself (see \cite{BL}). Hence,
  proving Theorem \ref{thm2.2} for each $s_i(\wedge^i\Pi_N)$ we will obtain that
 \begin{equation}\label{2.5}
 \lim_{N\to\infty}\frac 1N\ln s_1(\wedge^i\Pi_N)=\lim_{N\to\infty}\frac 1N\ln\prod_{j=1}^is_j
 (\Pi_N)=\sum_{j=1}^i\gam_j
 \end{equation}
 which yields (\ref{2.4}).

 The proof of Theorem \ref{thm2.2}, presented in Sections~\ref{sec3} and \ref{sec4}, is based on two
 main ingredients. The first one is a large deviations
 bound for products of random matrices which was first proved by Le Page under the additional
 contraction assumption. We rely on a version of this result from Theorem 14.19 in \cite{BQ} which does not require the contraction condition. In fact, the upper bound of large deviations from Theorem 6.2 on p.131 of \cite{BL} suffices for our purposes, as well.
 The second ingredient playing a decisive role in our proof of the lower bound below is the
 avalanche principle proved originally for two dimensional matrices in \cite{GS} and extended (in a strengthened form) to the multidimensional case in \cite{DK}. It is not difficult to see that the convergence in Theorem
  \ref{thm2.2} holds true also in mean which does not require large deviations estimates but only a subadditivity
  argument together with the avalanche principle.

 \subsection{Markov case} Next, we discuss the case when $\xi_0, \xi_1,\,\xi_2,...$ form a Markov chain on a
 Polish space $\cE$ (to conform with the standard notation, we start the indices from $0$), $F:\, \cE^\ell\to SL_d(\bbR)$ is a Borel measurable matrix function
 and $X_n=F(\xi_{q_1(n)},\xi_{q_2(n)},...,\xi_{q_\ell(n)})$ with $q_i(n),\, i=1,...,\ell$ satisfying
 Assumption \ref{ass2.1}(iii).  Let $P(n,x,\cdot),\,
 x\in \cE$ be the $n$-step transition probability of the Markov chain above, $P(x,\cdot)=
 P(1,x,\cdot)$ and assume that there exists
 a probability measure $\nu$ on $\cE$ such that for some $R,\rho>0$, all $n\geq 1$ and any
 bounded Borel function $f$ on $\cE$,
 \begin{equation}\label{2.6}
 \sup_{x\in \cE}|\int P(n,x,dy)f(y)-\int fd\nu|\leq Re^{-\rho n}\sup_{x\in \cE}|f(x)|.
 \end{equation}
 This assumption will be satisfied for an aperiodic Markov chain if, for instance, a version of the
 Doeblin condition holds true (see, for instance, \cite{Br}, Section 21.23).
 It follows that $\nu$ is the unique invariant measure of this Markov chain, i.e. the only measure $\nu$
 satisfying $\int d\nu(x)P(x,\Gam)=\nu(\Gam)$ for any Borel set $\Gam\subset \cE$, and so $\nu$ is
 ergodic. Taking $\nu$ as the initial distribution of the Markov chain, i.e. as the distribution of $\xi_0$,
  makes it a stationary ergodic process. Still, the condition (\ref{2.6}) will enable us to obtain stronger
  results for the Markov chain starting at any initial point $x\in\cE$.

 Let $\{\xi_n^{(i)},\, n\geq 0\},\, i=1,...,\ell$ be $\ell$ independent copies of the
  Markov chain $\{\xi_n,\, n\geq 0\}$ which produces an $\ell$-component Markov chain $\Xi_n=(\xi^{(1)}_{q_1(n)},\xi^{(2)}_{q_2(n)},...,\xi^{(\ell)}_{q_\ell( n)}),\, n\geq 0$ with the
  transition probabilities $P_\Xi(\bar x,\Gam_1\times\Gam_2\times\cdots\times\Gam_\ell)=
  \prod_{i=1}^\ell P(x_i,\Gam_i)$ where $\bar x=(x_1,...,x_\ell)$. Set $Y_n=F(\xi^{(1)}_{q_1(n)},
  \xi^{(2)}_{q_2(n)},...,\xi_{q_\ell( n)}^{(\ell)}),\, n\geq 0$ and assume that for some $\al>0$,
 \begin{equation}\label{2.7}
 \sup_{\bar x\in\cE^\ell}E_{\bar x}\| Y_1\|^\al<\infty
 \end{equation}
 where $E_{\bar x},\,\bar x=(x_1,...,x_\ell)$ is the expectation with respect to the
 probability $P_{\bar x}$ of the Markov chain $\Xi_n,\, n\geq 0$ starting at $\bar x$.

 Set $H_n=Y_n\cdots Y_2Y_1$. It follows from \cite{Bo2} (see also Section \ref{sec5}) that the limits
 \begin{equation}\label{2.8}
 \lim_{N\to\infty}\frac 1N\ln s_i(H_N)=\gam_i,\,\, i=1,...,d
 \end{equation}
 exist $P_{\bar x}$-almost surely (a.s.) for each $\bar x\in\cE^\ell$ where, again, $s_i(g)$ is the $i$-th
  singular value of a matrix $g$. Viewing (\ref{2.8}) as a definition of $\gam_i$'s we
 assume also that, for some $1 < k \leq d$,
 \begin{equation}
 \label{eq:sim} \gam_1>\gam_2 > \cdots > \gamma_k~;
 \end{equation}
 sufficient conditions for this
 can be found in \cite{Bo1} and \cite{Vi}. In addition, following \cite{Bo2} we assume
 quasi-irreducibility which means that the subspaces
 \[
 V(\bar x)=\{ u\in\bbR^d:\, \lim_{N\to\infty}\frac 1N\ln\| H_Nu\|\leq\gam_2\quad P_{\bar x}-
 \mbox{a.s.}\}
 \]
 are trivial for almost all $\bar x=(x_1,...,x_\ell)$ with respect to the product measure
 $\bar\nu$. Denote by $P_x$ the path space probability of the Markov chain $\xi_n,\, n\geq 0$
 provided that $\xi_0=x$.

 \begin{theorem}\label{thm2.3} Assume the above conditions (\ref{2.6}), (\ref{2.7}), (\ref{2.8}), (\ref{eq:sim}) and quasi-irreducibility. The singular values $s_i(\Pi_N),\,
 i=1,...,k-1$ of $\Pi_N=X_N\cdots X_2X_1$ satisfy
 \begin{equation}\label{2.9}
 \lim_{N\to\infty}\frac 1N\ln s_i(\Pi_N)=\gam_i\quad P_x-\mbox{a.s.}
 \end{equation}
 for each $x\in\cE$.
 \end{theorem}

 The proof of this result will be given in Section \ref{sec5} relying on the large deviations
 theorem for products of Markov dependent random matrices from \cite{Bo2} and an additional
 argument enabling us to compare large deviations estimates for the products $H_m$ and
 for $\Pi_{n+m}\Pi^{-1}_n$ in spite of the fact that the latter is not a product of Markov
 dependent random matrices.

 \section{Upper bound}\label{sec3}\setcounter{equation}{0}

 There are two cases in the proof of Theorem \ref{thm2.2}: $\gam_1=0$ and $\gam_1>0$.
 The first case requires only the upper bound since $\ln\| A\|\geq 0$ for any
 $A\in SL_d(\bbR)$. The second case will require both a lower and an upper bound
 so we will start with the latter which will serve in both cases. In fact, by Furstenberg's
 theorem (see Theorem 6.3 on p.66 in \cite{BL}) under the strong irreducibility condition $\gam_1=0$
 if and only if $G_\mu$ is contained in a compact subgroup; then each $X_n$ belongs to this subgroup
 too and Theorem \ref{thm2.2} follows in this case directly.

 It follows from the large deviations theorem for products of i.i.d.\ random matrices
 (see \cite[p.131, Theorem 6.2]{BL} and \cite[Theorem 14.19]{BQ}) that for any $\ve>0$ there exists $\ka(\ve)>0$
 and $n_1(\ve)\geq 1$ such that
 \begin{equation}\label{3.1}
 P\{\frac 1n\ln\| Y_n\cdots Y_2Y_1\|>\gam_1+\ve\}\leq e^{-\ka(\ve)n}
 \end{equation}
 for all $n\geq n_1(\ve)$. Without loss of generality $\kappa(\epsilon) < 1$. Fix $\ve>0$ and set $r(n)=r_\ve(n)=[\frac 2{\ka(\ve)}\ln n]$.
 Observe that if $r(n)\geq 1$ and $q_i(n+r(n))\leq q_{i+1}(n)$ for $i=1,...,\ell-1$
 then $X_n,X_{n+1},...,X_{n+r(n)-1}$ is an i.i.d.\ tuple having the same distribution
 as $(Y_1,Y_2,...,Y_{r(n)})$. Set $n_2(\ve)=\min\{ m\geq n_0(\frac 2{\ka(\ve)}):\, r(m)
 \geq n_1(\ve)\}$ where $n_0$ comes from Assumption \ref{ass2.1}(iii). Then for
 all $n\geq n_2(\ve)$,
 \begin{equation}\label{3.2}
 P\{\frac 1{r(n)}\ln\| X_{n+r(n)-1}\cdots X_{n+1}X_n\|>\gam_1+\ve\}\leq
 e^{-\ka(\ve)r(n)}.
 \end{equation}
 This together with (\ref{2.2}) and the Borel-Cantelli lemma yields existence of a
 finite with probability one random variable $M_{1}(\ve) = M_1(\ve,\om)$ such that for any $n\geq M_{1}(\ve)$,
 \begin{equation}\label{3.3}
 \ln\| X_{n+r(n)-1}\cdots X_{n+1}X_n\|\leq r(n)(\gam_1+\ve).
 \end{equation}

 Set $m_1=n_2(\ve)$ and recursively $m_{i+1}=
 m_i+r(m_i)$. Then $m_1<m_2<m_3<...$ and $m_n\to\infty$ as $n\to\infty$. Hence,
 $X_{m_i},X_{m_i+1},...,X_{m_{i+1}-1}$ is a tuple of i.i.d.\ random matrices
 for each $i\geq 1$. In particular, when $m_i\geq M_{1}(\ve)$ we have by (\ref{3.3}) that
 \begin{equation}\label{3.4}
 \ln\| X_{m_{i+1}-1}\cdots X_{m_i+1}X_{m_i}\|\leq (m_{i+1}-m_i)(\gam+\ve).
 \end{equation}
 By the submultiplicative property of the Euclidean operator matrix norm,
 \begin{equation}\label{3.5}
 \begin{split}
 &\ln\| X_N\cdots X_2X_1\|\leq \sum_{N\geq j\geq m_{k(N)}}\ln\| X_{j}\|\\
 & +\sum_{i:\, M_1(\ve)\leq m_i<k(N)}\ln\| X_{m_{i+1}-1}\cdots X_{m_i+1}X_{m_i}\|\\
 &+\sum_{1\leq j\leq\max(n_2(\ve),M_1(\ve)+r(M_1(\ve))}\ln\| X_j\|
 \end{split}
 \end{equation}
 where $k(N)=k_\ve(N)=\max\{ i:\, m_i\leq N\}$.
Since the last sum is a fixed random variable (depending on $\ve$) which is finite
with probability one then
\begin{equation}\label{3.6}
\lim_{N\to\infty}\frac 1N\sum_{1\leq j\leq\max(n_2(\ve),M_1(\ve)+r(M_1(\ve))}\ln\| X_j\|=0\quad\mbox{almost surely}.
\end{equation}

Next, we observe that by the Chebyshev inequality
\begin{equation}\label{3.7}
P\{\ln\| X_n\|\geq\frac 2\al\ln n\}=P\{\| X_1\|^\al\geq n^2\}\leq Dn^{-2}
\end{equation}
where $D=E\| X_1\|^\al<\infty$ by (\ref{2.1}). Hence, by the Borel-Cantelli lemma there exists a finite with
probability one random variable $M_2=M_2(\om)$ such that for all $n\geq M_2$,
\begin{equation}\label{3.8}
\ln\| X_n\|<\frac 2\al\ln n.
\end{equation}
Observe that
\begin{equation}\label{3.9}
N-m_{k(N)}<r(m_{k(N)})\leq r(N)=[\frac 2{\ka(\ve)}\ln N],
\end{equation}
and so, in particular, $m_{k(N)}\to\infty$ as $N\to\infty$. Thus, it suffices to estimate
the first expression in the right hand side of (\ref{3.5}) on the events $\Gam_N=
\{\om:\, M_2(\om)\leq m_{k(N)}\}$. By (\ref{3.8}) and (\ref{3.9}) on the event $\Gam_N$,
\begin{equation}\label{3.10}
\limsup_{N\to\infty}\frac 1N\sum_{N\geq j\geq m_{k(N)}}\ln\| X_{j}\|
\leq\limsup_{N\to\infty}\frac 2{N\al}r(N)\ln n=0.
\end{equation}

Finally, collecting (\ref{3.4})--(\ref{3.6}) and (\ref{3.8})--(\ref{3.10})
 we see that with probability one
\begin{equation*}
\limsup_{N\to\infty}\frac 1N\ln\| X_N\cdots X_2X_1\|\leq\gam_1+\ve.
\end{equation*}
Since $\ve>0$ is arbitrary we obtain the required upper bound
\begin{equation}\label{3.11}
\limsup_{N\to\infty}\frac 1N\ln\| \Pi_N \| = \limsup_{N\to\infty}\frac 1N\ln\| X_N\cdots X_2X_1\|\leq\gam_1
\end{equation}
with probability one. If $\gam_1=0$ this already implies (\ref{2.4}) while in the case
$\gam_1>0$ we shall  also need the corresponding lower bound.

\section{Lower bound}\label{sec4}\setcounter{equation}{0}

First, observe that without loss of generality we can assume here that $\gam_1>\gam_2$
where the $\gam_i$'s were defined in (\ref{2.3}). Indeed, either $\gam_1=\gam_2=
...=\gam_d$ and then $\gam_i=0$ for all $i$'s since all the matrices here have  determinant
equal one, or $\gam_1=...=\gam_k>\gam_{k+1}\geq...\geq\gam_d$ for some $1\leq k< d$. Then
we can prove the result for the first singular value of the $k$-th exterior power
$\wedge^k\Pi_N$ of $\Pi_N$ obtaining that with probability one,
\[
\lim_{N\to\infty}\frac 1N\ln s_1(\wedge^k\Pi_N) =\lim_{N\to\infty}\frac 1N
\sum_{i=1}^k\ln s_i(\Pi_N)=k\gam_1.
\]
Since $s_1(\Pi_N)\geq s_2(\Pi_N)\geq...\geq s_d(\Pi_N)$ and $\lim_{N\to\infty}\frac 1N
s_1(\Pi_N)\leq\gam_1$ with probability one by the upper bound we obtain that, in fact,
the last inequality is the equality. Thus, we obtain Theorem \ref{thm2.2} for
$s_1(\Pi_N)$ which is sufficient for its full statement as explained in Section \ref{sec2}.

Hence, we can and will assume here that $\gam_1>\gam_2,\,\gam_1>0$ and start with another bound
of large deviations for products of i.i.d.\ random matrices (see \cite{BQ}) which in the
same notation as in Section \ref{sec3} says that for any $\ve>0$ there exists $\ka(\ve)>0$ and
 $n_1(\ve)\geq 1$ such that
 \begin{equation}\label{4.1}
 P\{\frac 1n\ln\| Y_n\cdots Y_2Y_1\|<\gam_1-\ve\}\leq e^{-\ka(\ve)n}
 \end{equation}
 for all $n\geq n_1(\ve)$. Let $r(n)$ and $n_2(\ve)$ be the same as in Section \ref{sec3}.
 Then, for all $n\geq n_2(\ve)$ we obtain
 \begin{equation}\label{4.2}
 P\{\frac 1{r(n)}\ln\| X_{n+r(n)-1}\cdots X_{n+1}X_n\|<\gam_1-\ve\}\leq e^{-\ka(\ve)r(n)}.
 \end{equation}

 Since there exists no inequality similar to (\ref{3.5}) to employ for a proof of
 the lower bound we will need a more advanced argument in order to make use of the
 splitting of the product $X_N\cdots X_2X_1$ into appropriate products of i.i.d.\
 matrices. Namely, we will rely on the avalanche principle which appears for products
 of multidimensional matrices in \cite{DK}.
 Following \cite{DK} for each $g\in GL_d(\bbR)$ we set
 \[
 gr(g)=\frac {s_1(g)}{s_2(g)}
 \]
 which is called the gap of $g\in GL_d(\bbR)$. Now we have (see \cite[\S 2.4]{DK}),
 \begin{theorem}\label{thm4.1} (Avalanche Principle).
There exist universal constants $c,C>0$ such that whenever $a\geq cb>c$ and
$g_j\in GL_d(\bbR),\, j=1,...,l$ satisfy
\begin{enumerate}
\item[(i)]$gr(g_j)\geq a,\, j=1,...,l$ and
\item[(ii)]$\ln\| g_{j+1}g_j\|-\ln\| g_{j+1}\|-\ln\| g_j\|\geq -\frac 12\ln b$
\end{enumerate}
then
\begin{equation}\label{4.3}
\ln\| g_l\cdots g_2g_1\| +\sum_{j=2}^{l-1}\ln\| g_j\|\geq\sum_{j=1}^{l-1}
\ln\| g_{j+1}g_j\|-Cl\frac ba,\,\, j=1,...,l-1.
\end{equation}
\end{theorem}

Observe that from (ii) and (\ref{4.3}) we obtain
\begin{equation}\label{4.4}
\ln\| g_l\cdots g_2g_1\|\geq \sum_{j=1}^{l}\ln\| g_j\|-\frac 12l\ln b-Cl\frac ba.
\end{equation}
Let us take
\begin{equation}\label{4.5}
g_j=X_{m_{j+1}-1}\cdots X_{m_{j+1}}X_{m_j}
\end{equation}
where $m_1<m_2<...<m_{k(N)}$ are  as in Section \ref{sec3}. This together with (\ref{4.4})
will yield the required lower bound of the form
\begin{equation}\label{4.6}
\liminf_{N\to\infty}\frac 1N\ln\| X_N\cdots X_2X_1\|\geq\gam_1-\ve
\end{equation}
provided that we can obtain appropriate bounds on parameters $a$ and $b$ in the avalanche
principle above.

Now, (\ref{4.2}) together with the definition of $r(n)=r_\ve(n)$
and the Borel-Cantelli lemma yield that there exists
a finite with probability one random variable $M_{1}(\ve)$ such that for any $n\geq M_{1}(\ve)$,
\begin{equation}\label{4.7}
\ln\| X_{n+r(n)-1}\cdots X_{n+1}X_n\|\geq r(n)(\gam_1-\ve).
\end{equation}
In particular, for each $i<k(N)$ such that $m_i\geq M_{1}(\ve)$ we have
\begin{equation}\label{4.8}
\ln\| X_{m_{i+1}-1}\cdots X_{m_i+1}X_{m_i}\|\geq (m_{i+1}-m_i)(\gam_1-\ve).
\end{equation}

Next, set $j_{N}=\min\{ j:\, m_j\geq \sqrt N\}$. By the submultiplicative property of
the Euclidean matrix norm,
\begin{equation}\begin{split}\label{4.9}
&\ln\| X_N\cdots X_2X_1\| \geq\ln\| X_{m_{k(N)}-1}\cdots X_{m_{j_{N}}+1}
X_{m_{j_{N}}}\| \\
&\qquad -\ln\|(X_N\cdots X_{m_{k(N)+1}}X_{m_{k(N)}})^{-1}\|-\ln\|(X_{m_{j_{N}}-1}
\cdots X_2X_1)^{-1}\|.
\end{split}\end{equation}
As explained in Section 2 the condition (\ref{2.1}) implies also that $D'=E\| X_1^{-1}\|^{\al'}<\infty$
where $\al'=\frac \al{d-1}$. Thus, in the same way as in (\ref{3.7}) we have
\begin{equation}\label{4.10}
P\{\ln\| X_n^{-1}\|\geq\frac 2{\al'}\ln n\}\leq D'n^{-2},
\end{equation}
and so by the Borel-Cantelli lemma there exists a finite with
probability one random variable $M'_2=M'_2(\om)$ such that for all $n\geq M'_2$,
\begin{equation}\label{4.10+}
\ln\| X_n^{-1}\|<\frac 2{\al'}\ln n.
\end{equation}
Thus, similarly to (\ref{3.10}) we obtain that on the event $\Gam'_N=\{\om:\, M'_2(\om)\leq m_{k(N)}\}$,
\begin{equation}\label{4.11}
\limsup_{N\to\infty}\frac 1N\ln\| (X_N\cdots X_{m_{k(N)+1}}X_{m_{k(N)}})^{-1}\|
\leq\limsup_{N\to\infty}\frac 1N\sum_{N\geq j\geq m_{k(N)}}\ln\| X_{j}^{-1}\|=0.
\end{equation}
On the other hand,
\begin{equation}\label{4.11+}\begin{split}
&\limsup_{N\to\infty}\frac 1N\ln\|(X_{m_{j_{N}-1}}\cdots X_2X_1)^{-1}\|\\
&\leq\limsup_{N\to\infty}\frac 1N\sum_{1\leq j< M'_2}\| X_j^{-1}\|+
\limsup_{N\to\infty}\frac 1N\sum_{M'_2\leq j\leq m_{j_{N}}}\| X_j^{-1}\|=0.
\end{split}\end{equation}
Indeed, the first limit in the right hand side of (\ref{4.11+}) is zero since the sum there
is a fixed random variable which is finite with probability one. The second limit there is
zero in view of (\ref{4.11}) and the estimate $m_{j_{N}}\leq \sqrt N+r([\sqrt N])$.

Applying the avalanche principle we will show that in the above case with probability one
\begin{eqnarray}\label{4.12}
&\liminf_{N\to\infty}\frac 1N\ln\| X_{m_{k(N)}-1}\cdots
X_{m_{j_{N}}+1}X_{m_{j_{N}}}\|\\
&=\liminf_{N\to\infty}\frac 1N\ln\| g_{k(N)-1}\cdots g_{j_{N}+1}g_{j_{N}}\|
\geq \gam_1-7\ve. \nonumber
\end{eqnarray}
First, we estimate the avalanche principle parameters $a=a(\ve,N)$ and
$b=b(\ve,N)$ which will depend on $\ve$ and $N$. Set $g(n)=X_{n+r(n)-1}
\cdots X_{n+1}X_n$ so that $g_j = g(m_j)$, and let $s_1(g(n))\geq s_2(g(n))\geq...\geq s_d(g(n))>0$ be the
singular values of $g(n)$. The second exterior power $\wedge^2g(n)$ of $g(n)$ acting on
the second exterior power $\wedge^2\bbR^d$ of $\bbR^d$ has the biggest singular value
equal to $s_1(g(n))s_2(g(n))$. Hence
\begin{equation}\label{4.15}
gr(g(n))=\frac {s_1(g(n))}{s_2(g(n))}=\frac {s^2_1(g(n))}{s_1(g(n))s_2(g(n))}=
\frac {\| g(n)\|^2}{\|\wedge^2g(n)\|}
\end{equation}
where $\|\cdot\|$ is the Euclidean operator norm.

Now, set $H_n=Y_n\cdots Y_2Y_1$ with $Y_1,Y_2,...$ introduced in Section \ref{sec2}.
Under our conditions with probability one
\begin{equation}\label{4.16}
\lim_{n\to\infty}\frac 1n\ln\| H_n\|=\gam_1\,\,\,\mbox{and}\,\,\,\lim_{n\to\infty}
\frac 1n\ln\|\wedge^2H_n\|=\gam_1+\gam_2
\end{equation}
and, recall that $\gam_2<\gam_1$. Applying the large deviations bounds
to $\| H_n\|$ and to $\|\wedge^2H_n\|$ we obtain that for any $\ve>0$ there exists
$\ka(\ve)>0$ (which could be different from before but we denote it by the same letter)
and $n_3(\ve)\geq 1$ such that
\begin{equation}\label{4.17}
P\{\frac 1n\ln\|\wedge^2H_n\|>\gam_1+\gam_2+\ve\}\leq e^{-\ka(\ve)n}
\end{equation}
for all $n\geq n_3(\ve)$. Hence, if $r(n)\geq n_3(\ve)$ and $n\geq n_0(\frac 2{\ka(\ve)})$
then
\begin{equation}\label{4.18}
P\{\frac 1{r(n)}\ln\|\wedge^2g(n)\|>\gam_1+\gam_2+\ve\}\leq e^{-\ka(\ve)r(n)}.
\end{equation}
This together with (\ref{4.2}) and (\ref{4.15}) yields that
\begin{equation}\label{4.19}
P\{ gr(g(n))<e^{(\gam_1-\gam_2-2\ve)r(n)}\}\leq 2e^{-\ka(\ve)r(n)}.
\end{equation}
Taking into account that $r(n)=[\frac 2{\ka(\ve)}\ln n]$ we conclude from
(\ref{4.19}) and the Borel-Cantelli lemma that there exists a finite with probability
one random variable $M_{3}(\ve)$ such that for any $n\geq M_{3}(\ve)$,
\begin{equation}\label{4.20}
gr(g(n))\geq e^{(\gam_1-\gam_2-2\ve)r(n)}.
\end{equation}

Next, we use that by our choice of $r(n)$ there exists $n_4(\ve)\geq 1$
such that if $n\geq n_4(\ve)$ then
$X_n,X_{n+1},...,X_{n+r(n)+r(n+r(n+r(n)))-1}$ is an i.i.d.\ tuple
having the same distribution as $Y_1,Y_2,...,Y_{r(n)+r(n+r(n+r(n)))}$.
Thus, similarly to the above, relying on the large deviations bound
(\ref{4.1}) together with the Borel-Cantelli lemma we conclude that
there exists a finite with probability one random variable $M_{4}(\ve)$
such that for any $n\geq M_{4}(\ve)$,
\begin{eqnarray}\label{4.21}
&\| X_{n+r(n)+r(n+r(n+r(n)))-1}\cdots X_nX_{n+1}\|\\
&=\| g(n+r(n))g(n)\|\geq e^{(\gam_1-\ve)(r(n)+r(n+r(n)))}.\nonumber
\end{eqnarray}
Applying the large deviations estimate (\ref{3.3}) to $\| g(n)\|$ and
to $\| g(n+r(n))\|$ together with the Borel-Cantelli lemma we obtain
that there exists a finite with probability one random variable $M_{5}(\ve)$
 such that for any $n\geq M_{5}(\ve)$,
\begin{equation}\label{4.22}
\| g(n)\|\leq e^{(\gam_1+\ve)r(n)}\,\,\,\mbox{and}\,\,\,\| g(n+r(n))\|
\leq e^{(\gam_1+\ve)r(n+r(n))}.
\end{equation}
Let $n_5(\ve)$ be such that $\frac 2{\ka(\ve)}\ln(1+\frac {r_\ve(n)}n)\leq 1$
for any $n\geq n_5(\ve)$.
Then, by (\ref{2.2}), (\ref{4.21}) and (\ref{4.22}) for any $n\geq\max( n_5(\ve), M_{5}(\ve))$,
\begin{equation}\label{4.23}
\frac {\| g(n+r(n))g(n)\|}{\| g(n)\|\| g(n+r(n))\|}\geq
e^{-3\ve(r(n)+r(n+r(n)))}\geq e^{-6\ve(r(n)+1)}.
\end{equation}

Observe that for $n\geq \sqrt N$ the numbers $k(N)=\max\{ i:\, m_i<N\}$ and $j_{N} = \min \{ j \, : \, m_j \geq \sqrt N\}$ satisfy
\begin{equation}\label{4.24}
k(N)-j_{N}\leq\frac N{r([\sqrt N])}=\frac N{[\frac {1}{\ka(\ve)}\ln N]}.
\end{equation}
When $n\geq \sqrt N$ then (\ref{4.20}) and (\ref{4.23}) hold true for
\[
\om\in\Om_{\ve,N}=\{\om:\,\max(M_{3}(\ve,\om),M_{4}(\ve,\om),M_{4}(\ve,\om)(\om),n_5(\ve))
\leq\sqrt N\}.
\]
Clearly, $\Om_{\ve,N}\uparrow\tilde\Om$ with $P(\tilde\Om)=1$. Thus we
can estimate the parameters of the avalanche principle for $\om\in
\Om_{\ve,N}$ and each fixed $N$ large enough and then let $N\to\infty$.

It follows from (\ref{4.20}), (\ref{4.23}) and (\ref{4.24}) that applying
the avalanche principle to $g_{j_{N}},\, g_{j_{N}+1},..., g_{k(N)-1}$
 we can take in (\ref{4.4}),
\begin{eqnarray*}
&l=l(N)=k(N)-j_{N},\, a=a(\ve,N)=e^{(\gam_1-\gam_2-2\ve)r([\sqrt N])}\\
&\mbox{and}\,\,\, b=b(\ve,N)=e^{6\ve (r(N)+1)}.
\end{eqnarray*}
Choosing $\ve$ much smaller than $\frac 18(\gam_1-\gam_2)$
we let $N\to\infty$ to obtain from (\ref{4.4}), (\ref{4.8}), (\ref{4.20}), (\ref{4.23})
and (\ref{4.24}) together with the avalanche principle that (\ref{4.12})
 holds true. These together with (\ref{4.9}), (\ref{4.11}) and (\ref{4.11+}) yield that
\begin{equation}\label{4.25}
\liminf_{N\to\infty}\frac 1N\| X_N\cdots X_2X_1\|\geq \gam_1-7\ve.
\end{equation}
Now we let $\ve\to 0$ and obtain
\[
\liminf_{N\to\infty}\frac 1N\ln\| X_N\cdots X_2X_1\|\geq\gam_1
\]
which together with (\ref{3.11}) yields (\ref{2.4}) and completes the proof
of Theorem \ref{thm2.2}.
\qed

\section{Products with Markov dependence}\label{sec5}\setcounter{equation}{0}

As in the case of Theorem \ref{thm2.2} it suffices to prove Theorem \ref{thm2.3} only
for the biggest singular value $s_1(\Pi_N)$. It is easy to see that the condition of the
form (\ref{2.6}) remains true also for the product Markov chain $\Xi_n,\, n\geq 0$. Hence,
 it follows by the large deviations result of Theorem 4.3
in \cite{Bo2} applied to the products $H_N=Y_N\cdots Y_2Y_1$ of Markov dependent random
matrices that for any $\ve>0$ there exists $\ka(\ve)>0$ and $n(\ve)\geq 1$ such that
\begin{equation}\label{5.1}
P_{\bar x}\{|\frac 1n\ln\| H_n\|-\gam_1|>\ve\}\leq e^{-\ka(\ve)n}
\end{equation}
for any $n\geq n(\ve)$ and $\bar x=(x_1,...,x_\ell)\in E^\ell$ where, recall, $P_{\bar x}$
is the probability conditioned on $\Xi_0=\bar x$. Note that (\ref{5.1}) together with the
Borel--Cantelli lemma yields (\ref{2.7}).

Next, we observe that (\ref{2.6}) implies $\phi$-mixing of the Markov chain $\xi_0,\xi_1,
\xi_2,...$ with the $\phi$-dependence coefficient satisfying $\phi(n)\leq 2Re^{-\rho n}$ (see \cite{Br}).
This seems to be well known (see p.p.365-366 in \cite{IL} for the case of finite Markov chains
and Theorem 21.1 in \cite{Br} for a general stationary Markov chain) but we claim this for
each probability $P_x,\, x\in\cE$, and so
for readers' convenience we will elaborate this here. For any $0\leq m\leq n$ let $\cF_{mn}$
be the $\sig$-algebra generated by $\xi_m,\xi_{m+1},...,\xi_n$. Then the $\phi_x$-dependence
coefficient for $x\in\cE$ is defined by
\begin{equation}\label{5.2}
\phi_x(n)=\sup_{m\geq 0}\big\{\big\vert\frac {P_x(\Gam\cap\Del)}{P_x(\Gam)}-P_x(\Del)
\big\vert :\, \Gam\in\cF_{0m},\,\Del\in\cF_{m+n,\infty},\, P_x(\Gam)>0\big\}
\end{equation}
where, recall, $P_x$ is the probability corresponding to the initial condition
 $\xi_0=x$. In order to show that
\begin{equation}\label{5.3}
\phi_x(n)\leq 2Re^{-\rho n}
\end{equation}
when (\ref{2.6}) holds true observe that it suffices to consider $\Gam$ and $\Del$ of the
form
\[
\Gam=\bigcap_{i=1}^k\{\xi_{m_i}\in G_i\}\,\,\,\mbox{and}\,\,\,\Del=\bigcap_{j=k+1}^l
\{\xi_{m_j}\in G_j\}
\]
where $m_1<m_2<...<m_k<m_k+n\leq m_{k+1}<m_{k+2}<...<m_l$. Set
\begin{eqnarray*}
&f(y)=\bbI_{G_{k+1}}(y)\int_{G_{k+2}}P(m_{k+2}-m_{k+1},y,dy_1)\int_{G_{k+3}}P(m_{k+3}-m_{k+2},
y_1,dy_2)\\
&\cdots\int_{G_{l-1}}P(m_l-m_{l-1},y_{l-k-2},G_l),
\end{eqnarray*}
where $\bbI_G$ is the indicator of $G$, and observe that $P_\nu(\Del)=\int_Ef(x)d\nu(x)$.
Then
\begin{eqnarray*}
&P_x(\Gam\cap\Del)=\int_{G_1}P(m_1,x,dy_1)P(m_2-m_1,y_1,dy_2)\int_{G_2}P(m_3-m_2,y_2,dy_3)\int_{G_3}\\
&\cdots P(m_k-m_{k-1},y_{k-1},dy_k)\int_{G_k}P(m_{k+1}-m_k,y_k,dy_{k+1})\int_{G_{k+1}}
f(y_{k+1})\\
&=P_x(\Gam)P_x(\Del)+Q
\end{eqnarray*}
where by (\ref{2.6}),
\begin{eqnarray*}
&|Q|\leq P_x(\Gam)\sup_y\left|\int_\cE P(m_{k+1}-m_k,y,dz)f(z)-\int_\cE P(m_{k+1},x,dz)f(z)\right|\\
&\leq P_x(\Gam)\sup_y\left|\int_\cE P(m_{k+1}-m_k,y,dz)f(z)-\int f(z)d\nu(z)\right|\\
& +P_x(\Gam)\left|\int_\cE P(m_{k+1},x,dz)f(z)-\int f(z)d\nu(z)\right|\leq 2Re^{-\rho n}P_x(\Gam)
\end{eqnarray*}
yielding (\ref{5.3}).

Next, set $r(n)=r_\ve(n)=[\frac 2{\del(\ve)}\ln n]$, where $\del(\ve)=\min(\ka(\ve),\rho)$,
 and observe that for large $n$,
\begin{equation}\label{5.4}
q_{i+1}(n)\geq q_i(n+2r(n))\geq q_i(n+r(n))+r(n)\,\,\,\,\,\mbox{for all}\,\, i=1,...,\ell-1.
\end{equation}
Consider vectors $\bar x^{(i)}=(x_1^{(i)},x_2^{(i)},...,x^{(i)}_m),\, x_j^{(i)}\in\cE,\,
i=1,...,\ell$ and view products
\[
Q_m(\bar x^{(1)},\bar x^{(2)},...,\bar x^{(\ell)})=\prod_{j=1}^mF(x_j^{(1)},x_j^{(2)},...,
x_j^{(\ell)})
\]
as functions of vectors $\bar x^{(i)},\, i=1,...,\ell$. Introduce another function
\[
\vf_m(\bar x^{(1)},...,\bar x^{(\ell)})=\bbI_{\{ |\frac 1m\ln\| Q_m(\bar x^{(1)},...,
\bar x^{(\ell)})\|-\gam_1|>\ve\}}
\]
which takes on the values 0 and 1 only. We are going to plug in place of $\bar x^{(i)}$ in $\vf_m$
with $m=r(n)$ the vectors $\bar\xi^{(i)}=(\xi_{q_i(n)},\xi_{q_i(n+1)},...,\xi_{q_i(n+r(n)-1)})$
observing that $\bar\xi^{(i)}$ is $\cF_{q_i(n),q_i(n+r(n)-1)}$-measurable and that by (\ref{5.4})
there is a gap of at least $r(n)$ between the intervals $[q_i(n),\, q_i(n+r(n)-1)]$ for
different $i=1,...,\ell$.

To use the above observation we will need the following result which is a particular case
of Corollary 1.3.11 in \cite{HK} (see also Corollary 3.3 in \cite{Ha}).
\begin{lemma}\label{lem5.1}
Let $Z_i$ be $\wp_i$-dimensional $\cE^{\wp_i}$-valued random vectors with a distribution
$\mu_i,\, i=1,...,k$ defined on the same probability space $(\Om,\cF,P)$ and such that
$Z_i$ is $\cF_{m_in_i}$-measurable where $n_{i-1}<m_i\leq n_i<m_{i+1},\, i=1,...,k,\,
n_0=0,\, m_{k+1}=\infty$. Then for any bounded Borel function $h=h(x_1,...,x_k),\,
x_i\in\cE^{\wp_i}$,
\begin{eqnarray}\label{5.5}
&|Eh(Z_1,Z_2,...,Z_k)-\int h(x_1,x_2,...,x_k)d\mu_1(x_1)d\mu_2(x_2)...d\mu(x_k)|\\
&\leq 4\sup_{x_1,...,x_k}|h(x_1,...,x_k)|\sum_{i=2}^k\phi(m_i-n_{i-1})\nonumber
\end{eqnarray}
with the $\phi$-dependence coefficient defined in (\ref{5.2}). In particular, if
$Z_1^{(1)},Z_2^{(2)},...,Z_k^{(k)}$ are independent copies of $Z_1,Z_2,...,Z_k$,
respectively, then taking $h=\bbI_\Gam$ for a Borel set $\Gam\subset\cE^{\wp_1+\wp_2+
\cdots +\wp_k}$ it follows that
\begin{equation}\label{5.6}
|P\{ (Z_1,Z_2,...,Z_k)\in\Gam\}-P\{(Z_1^{(1)},Z_2^{(2)},...,Z_k^{(k)})\in\Gam\}|\leq
4\sum_{i=2}^k\phi(m_i-n_{i-1}).
\end{equation}
\end{lemma}

Now, applying (\ref{5.6}) to $Z_i=\bar\xi^{(i)},\, i=1,...,\ell$ and
\[
\Gam=\big\{(\bar x^{(1)},...,\bar x^{(\ell)}):\, |\frac 1{r(n)}\ln\| Q_{r(n)}(\bar x^{(1)},
...,\bar x^{(\ell)})\|-\gam_1 |>\ve\big\}
\]
and taking into account (\ref{5.1}), (\ref{5.3}) and (\ref{5.4}) we obtain that for each $x\in\cE$,
\begin{equation}\label{5.7}
P_x\big\{ |\frac 1{r(n)}\ln\| X_{n+r(n)-1}\cdots X_{n+1}X_n\|-\gam_1|>\ve\}\leq
e^{-\ka(\ve)r(n)}+8R\ell n^{-2}
\end{equation}
whenever $r(n)\geq n(\ve)$.

The remaining part of the proof of Theorem
\ref{thm2.3} proceeds in the same way as in the i.i.d.\ case of Theorem \ref{thm2.2}
except for the arguments leading to (\ref{3.10}), (\ref{4.11}) and (\ref{4.11+}). Namely,
we cannot use the Chebyshev inequality in order to obtain (\ref{3.7}) and (\ref{4.10})
since, in general, in the present situation $E_\nu\| X_n\|^\al$ and $E_\nu\| X_n^{-1}\|^{\al'}$ may be
not equal to $E_\nu\| X_1\|^\al$ and $E_\nu\| X_1^{-1}\|^{\al'}$, respectively, and the latter expectations
may be not equal to $E_{\bar\nu}\| Y_1\|^\al$ and $E_{\bar\nu}\| Y_1^{-1}\|^{\al'}$ where $E_\nu$ and $E_{\bar\nu}$ are the expectations corresponding to the path space probabilities $P_\nu$ and $P_{\bar\nu}$ of the
Markov chains $\xi_n$ and $\Xi_n$ having initial distributions $\nu$ and $\bar\nu=\nu\times\cdots\nu$,
 respectively. But applying Lemma \ref{lem5.1} in the same way as in (\ref{5.7})
we obtain by the Chebyshev inequality that
\begin{eqnarray*}
&P_x\{\ln\| X_n\|\geq\frac 2\al\ln n\}\leq P_{\bar x}\{\ln\| Y_n\|\geq\frac 2\al\ln n\}+8R\ell n^{-2}\\
&\leq n^{-2}E_{\bar x}\| Y_n\|^\al +8R\ell n^{-2}\leq n^{-2}(\sup_{\bar z}E_{\bar z}\| Y_1\|^\al+8R\ell)\quad
\forall\,\, x\in\cE
\end{eqnarray*}
since
\[
E_{\bar x}\| Y_n\|^\al=\int_\cE P_{\Xi}(n,\bar x,d\bar y)\| F(\bar y)\|^\al\leq\sup_{\bar z}
\int P_\Xi(\bar z,d\bar v)\| F(\bar v)\|^\al=\sup_{\bar z}E_{\bar z}\| Y_1\|^\al<\infty.
\]
Similarly, for any $x\in\cE$,
\begin{equation*}
P_x\{\ln\| X^{-1}_n\|\geq\frac 2{\al'}\ln n\}\leq (D'+8R\ell)n^{-2}
\end{equation*}
where $D'=E_{\bar x}\| Y^{-1}_1\|^{\al'}<\infty$ and $\al'=\frac \al{d-1}$. Now, the
corresponding versions of (\ref{3.10}), (\ref{4.11}) and (\ref{4.11+}) follow in the same way as in
Section \ref{sec3} and \ref{sec4} while the arguments related to the avalanche principle remain the same.
\qed
\begin{remark}\label{rem5.2}
Since Lemma \ref{lem5.1} is quite general the Markov dependence in the sequence
$\xi_n,\, n\geq 0$ is needed only to rely on large deviations result from \cite{Bo2},
and so our method will go through whenever large deviations estimates (actually, only
upper bounds in the form (\ref{5.1})) for products of
stationary sufficiently fast $\phi$-mixing sequences of random matrices become available.
\end{remark}

\bibliography{matz_nonarticles,matz_articles}
\bibliographystyle{alpha}

\end{document}